\documentclass{elsarticle}

\usepackage{amssymb}
\usepackage{amsmath}

\journal{J. Math. Anal. and Appl.}

\newtheorem{thm}{Theorem}
\newtheorem{prop}[thm]{Proposition}
\newtheorem{cor}[thm]{Corollary}
\newtheorem{lem}[thm]{Lemma}

\newtheorem{rem}[thm]{Remark}

\newproof{pf}{Proof}

\newcommand{\be}{\begin{equation}}
\newcommand{\ee}{\end{equation}}

\newcommand{\A}{\mathcal A}

\newcommand{\N}{{\mathbb{N}}}
\newcommand{\R}{{\mathbb{R}}}

\begin{document}

\begin{frontmatter}

\title{Rolewicz-type chaotic operators\tnoteref{t1}}
\tnotetext[t1]{This research was supported by the grant Cori 2013 of the University of Palermo. The second Author thanks the hospitality of the Department of Mathematics of University of Palermo}

\author[rvt]{D. ˜Bongiorno\corref{cor1}}

\ead{donatella.bongiorno@unipa.it}
\author[focal]{U.B. ˜Darji}

\ead{ubdarj01@louisville.edu}

\author[els]{L. ˜Di Piazza}

\ead{luisa.dipiazza@unipa.it}

\cortext[cor1]{Corresponding author}

\address[rvt]{Dipartimento di Energia, Ingegneria dell'informazione e Modelli matematici (DEIM), University of Palermo,  Palermo, Italy}

\address[focal]{Department of Mathematics, University of Louisville,  Louisville,  KY 40292, USA}

\address[els]{Dipartimento di Matematica e Informatica, University of Palermo, Via Archirafi 34, 90123 Palermo, Italy}

%\subjclass[2010]{Primary: 46G10, 28B05 Secondary: 15A03}
\maketitle

\begin{abstract} In this article we introduce a new class of Rolewicz-type operators in $\ell_p$, $1 \le p < \infty$. We exhibit a collection ${\mathcal F}$ of cardinality continuum of operators of this type which are chaotic and remain so under almost all finite linear combinations, provided that the linear combination has sufficiently large norm. As a corollary to our main result we also obtain that there exists a countable collection of such operators whose all finite linear combinations are chaotic provided that they have sufficiently large norm.
\end{abstract}

\end{frontmatter}
\date{}

\section*{Introduction}

Hypercyclic operators are generalizations of cyclic operators which have been studied in operator theory for many years. If $X$ is a Banach space and $T:X \rightarrow X$ is a bounded linear operator, we say that {\em $T$ is
hypercyclic} if there exists $x \in X$ such that $Orb(x,T)$, the orbit of $x$ under $T$, is dense in $X$. It is a well-known result of Rolewicz \cite{rolewicz} that no finite dimensional Banach space admits a hypercyclic operator. However, every infinite dimensional Banach space admits a hypercyclic operator.  This was shown independently by Ansari \cite{ansari}  and Bernal-Gonz\`alez \cite{bernal1}.

From the dynamical systems point of view, hypercyclic operators are closely connected to transitive maps. A continuous self-mapping $f$ of a metric space $X$ is {\em transitive} if for all non-empty open sets $U, V$
of $X$, there exists $n \ge 1$ such that $f^n(U) \cap V \neq \emptyset$. If $X$ is a separable metric space without isolated points, then $f$ is transitive if and only if $Orb(x, f)$ is dense in $X$ for some $x \in X$. Hence, in the case of a separable Banach space, a linear operator is hypercyclic if and only if it is transitive. The notion of chaos in the sense of Devaney \cite{devaney} also applies in the setting of Banach spaces. We say that a  bounded linear operator $T: X \rightarrow X$ is {\em chaotic in the sense of Devaney} if $T$ is transitive and the set of periodic points of $T$ is dense in $X$. In this article, when we say an operator is chaotic, we mean that the operator is chaotic in the sense of Devaney.
Whereas every infinite dimensional Banach space admits a hypercyclic operator, Bonet, Mart\'inez-Gim\'enez, and Peris
\cite{bonet1} showed that the separable hereditarily indecomposable Banach space constructed by Gowers and Maury \cite{gowers} admits no chaotic operators.

When one encounters an exotic or an unusual object, a natural question one asks is how many such objects are there. In topological setting, one asks if the collection of objects in questions forms a meager or comeager set. In every  separable infinite dimensional Hilbert space, the set of hypercyclic operators forms a nowhere dense set with respect to the norm topology  \cite[Thm 2.24]{bayartmatheron}. Hence, in some sense there are very few chaotic operators in infinite dimensional Hilbert space. It is easy to construct two chaotic operators whose sum is not chaotic. In fact, Grivaux \cite{grivaux} showed that every bounded linear operator on a separable  infinite dimensional Hilbert space is the sum of  two chaotic operators. Hence being chaotic is far from being stable under finite linear combinations. However, in recent years a new algebraic notion of largeness has been popularized and exploited
\cite{gurariy2}, \cite{aron1}, \cite{aron2}. The idea is to exhibit a large algebraic structure in a setting where no algebraic structure is apparent. If $X$ is a Banach space and $A \subseteq X$, we say that {\em $A$ is lineable}
if there is a vector space $V \subseteq X$ of dimension continuum such that every non-zero element of $V$
is in $A$. Recently, many sets of classical importance are shown to be lineable. For example, the set of nowhere differentiable functions is lineable \cite{gurariy1}. For other examples and survey on the subject refer to \cite{bdd} \cite{enflo}. Clearly, the set of chaotic operators cannot be lineable as no operator with norm less than one is hypercyclic. However, as we will see from the main result of this article, some type of algebraic largeness
can be inferred.

The notion of shift is a basic yet fundamental concept in dynamical systems. In symbolic dynamics they are of utmost importance and have been extensively studied. In 1969 Rolewicz \cite{rolewicz}  made an important observation. If one
considers the shift $T: \ell _2 \rightarrow \ell_2$, defined by $T(x_1, x_2, \dots) = (x_2, x_3, \dots)$, then
for all $\lambda >1$, $\lambda T$ is a transitive operator on $\ell_2$ whose set of periodic points is dense in $\ell_2$. Operators of these type are called Rolewicz operators. They have been generalized to
what are called weighted shift operators or weighted Rolewicz operators. For a suitable choice of sequence $\{w_n\}$, the shift operator with the weight sequence $\{w_n\}$ is defined by
\[ T(x_1, x_2, \dots) = (w_1x_2, w_2 x_3, \dots).
\]
Weighted shift operators have been well-studied. In this note, we view the shift operator in a different light. We think of the shift operator as the projection of the input vector on coordinates $\{2, 3, \dots\}$. In general, given an increasing function $f: \N \rightarrow \N$, we define $T_f: \ell_p \rightarrow \ell_p$, $ 1 \le p < \infty$, by
\[ T_f(x_1, x_2, \dots) =(x_{f(1)}, x_{f(2)}, \dots).
\]
Hence, the usual shift operator equals $T_f$, where $f(n) = n+1$. As in the Rolewicz operator, if $\lambda >1$,
then we have that $\lambda T_f$ is chaotic.  We call operators of this type {\em Rolewicz-type operators.}. We
study and exploit these operators. We prove that for suitable $f_1, \dots, f_t$,  $c_1, \dots, c_t$, and sufficiently large
$\lambda$, $\lambda \sum_{i=1}^tc_iT_{f_i}$ is chaotic. The precise formulation is given in Theorem~\ref{teore}.
As corollaries, we obtain the following results which show that set of chaotic operators contains large algebraic structure.

{\bf Corollary \ref{c1}}
{\it There exists  an infinite family ${\mathcal T}$ of chaotic operators such that for all $T_1,T_2, \dots , T_t\in{\mathcal T}$ and $c_1,c_2, \dots , c_t\in\R$ the operator $\lambda\sum_{i=1}^tc_iT_{f_i}$ is chaotic for sufficiently large $\lambda,$ provided that the following sum $\sum_{i=1}^tc_iT_{f_i}$ is not the zero operator.}

{\bf Corollary \ref{c2}}
{\it There exists a family ${\mathcal T}$ of cardinality continuum of chaotic operators such that for almost all $(c_1,c_2, \dots ,c_t)\in\R^t$ and $T_1, T_2, \dots , T_t\in \mathcal T$ the operator $\lambda\sum_{i=1}^tc_iT_{f_i}$ is chaotic for sufficiently large $\lambda.$}

We would like to point out that recently other notions of chaos have enjoyed attention in linear dynamics as well. These notions include Li-Yorke chaos and distributional chaos.  A characterization of Li-Yorke
chaos in Banach spaces in terms of irregular vectors was given in \cite{bermudez} and later extended to the setting of Fr\'echet spaces in \cite{bernardes1}.  As an easy consequence of these results one obtains that all hypercyclic operators are Li-Yorke chaotic. Hence, all operators constructed in our main results are Li-Yorke chaotic.
A characterization of distributional chaos in the setting of Fr\'echet spaces  in terms of distributionally irregular vectors was given in \cite{bernardes2}. We conjecture that all operators constructed in our main results are distributionally chaotic. Indeed, we conjecture
that if $f_1, \dots f_t$ are strictly increasing pairwise almost disjoint functions  and $c_1, \dots, c_t$ are real numbers, then $T= \sum_{i=1}^t c_i T_{f_i}$ being chaotic (in the sense of Devaney) implies that $T$ is
distributionally chaotic. Indeed, for weighted backward shifts it was shown in \cite{martinez1} that Devaney chaos implies distributional chaos. We conjecture that using the characterization of distributional chaos given in \cite{bernardes2}, one may be able to prove a large class of Devaney chaotic operators, which include weighted shifts as well as Rolewicz-type operators defined in this article, are distributionally chaotic.

\section{Notations and preliminaries}
In what follows $X$ is a separable  infinite dimensional Banach space. Let $T:X\to X$ be a linear bounded   operator. For $x \in X$ the {\it orbit} of
$x$ under $T$ is ${\rm Orb}(T,x)= \{x, T(x), T^2(x),...\}$ where $T^n=T \circ T \circ ...\circ T$ is the $n^{th}$ iterate of $T$ obtained by composing
$T$ with itself n times. A point $x \in X$ is a {\it periodic point} of $T$ if $T^n(x)=x$ for some $n \geq 1$. Also $T$ is {\it transitive} if for any two non-empty open sets $U, \ V$ in $X$, there exists an integer $n \geq 1$ such that $T^n(U) \cap  V \neq \emptyset$.

\bigskip
{\bf Devaney's Definition of Chaos} (\cite{devaney}) We say that a linear   bounded operator $T:X\to X$ is {\it chaotic} if
\begin{enumerate}
\item[(i)]   $T$ is transitive;
\item[(ii)] the periodic points of $T$ are dense in $X$.
\end{enumerate}
It is well known that, in a complete metric space with no isolated points, being transitive is equivalent (via the Baire Category Theorem) to having a
point with dense orbit, which it is equivalent to having a dense $G_{\delta}$ set of points each of which has a dense orbit.\\

We refer the reader to texts \cite{bayartmatheron} and \cite{grosse} for further information on linear chaos.\\

From now on we will consider linear bounded operators in $l_p$,  for some fixed $p$ where $1\leq p<\infty$. If $ \vec{x}=(x_1,x_2,\dots ) \in l_p$ we  denote by $\|\vec{x}\|_p:= (\sum_{i=1}^{\infty} |x_i|^p)^{\frac{1}{p}}$ and  by $\|\vec{x}\|_{\infty}=\sup\{|x_i|: \ i=1,2,\dots\}, $ respectively its $l_p$  and sup norm.\\

 We start with some preliminary facts and notations.\\

Let $\N= \{1,2,\dots\}$.  We say that a map $f:\N\to\N$ is  {\it  increasing} if   $f(i) <f(j)$ whenever $i<j$.
Let  $f:\N\to\N $ and $g:\N\to\N $ be given. We say that $f$ and $g$ are {\it almost disjoint} if the set
$\{k:f(k)=g(k)\}$ is finite and  $f(k)\neq g(k')$ if $k\neq k'$.\\

For a given increasing function  $f:\N\to\N$, denote by
 $T_{f}:l_p\to l_p$ the linear bounded operator   defined by
 \[ T_{f}(x_1,x_2,\dots ) := (x_{f(1)},x_{f(2)}, \dots ). \]

Let $f_1, \ f_2,  \dots, f_t$ be increasing functions  and let $c_1, \ c_2, \dots, c_t \in
\R $. Then for each $\sigma \in \{1, \ 2, \dots, \
t\}^n$ we define $$f_{\sigma}:=f_{\sigma(n)}\circ f_{\sigma(n-1)}
\circ \dots\circ f_{\sigma(1)}$$ and
$$c_{\sigma}:=c_{\sigma(n)}\cdot  c_{\sigma(n-1)}\cdot \dots \cdot
c_{\sigma(1)}.$$

If $\sigma, \tau $ are in $\{1, \ 2, \dots, \ t\}^n$ and
$\{1, \ 2, \dots, \ t\}^m$, respectively, then $\sigma \tau$ is
just the concatenation of $\sigma$ followed by $\tau$, and, in general, if $S \subseteq \{1, \ 2, \dots, \ t\}^n$ then $\sigma S = \{\sigma \tau : \tau \in S \}$. We use $|\sigma|$ to denote the length of $\sigma$. In
particular $f_{\sigma}f_{\tau}=f_{\sigma \tau}$ and  $c_{\sigma} c_{\tau}=c_{\sigma \tau}$. For notational convenience we let $f_{\emptyset}= id$ and $c_{\emptyset}= 1$. If $S \subset \{1,2,\dots , t\}^r,$ then we let $c(S):=\sum_{\sigma \in S}c_{\sigma}.$ Now we define a collection of equivalence relations. Fix
$r, i \in \N$. Let $\sigma , \tau\in \{1,2,\dots , t\}^r$.

We say that $\sigma \sim_i \tau$ if and only if $f_{\sigma}(i)=f_{\tau}(i).$
We note that $\sim_i$ is an equivalence relation on $\{1,2,\cdots , t\}^r$.
If $\sigma \in \{1,2,\dots , t\}^r$, then
\[[\sigma]_i =\{ \tau\in\{ 1,2,\dots ,t\}^r:\sigma\sim_i\tau\}. \]
We note that in general there is no ambiguity in writing $[\sigma]_i$ as $r = |\sigma|$.
We let $\left| [\sigma]_i \right|$ the number of elements of $[\sigma]_i$ and
\[\A (r,i) = \{ [\sigma]_i: \sigma \in \{1, \dots, t \}^r\},
\]
and $\sharp (r,i)$ the cardinality of $\A(r,i)$.

We say that the constants $c_1,c_2, \dots , c_t \in\R$ satisfy the
{\em non-zero condition at level $m$, $m \in \N$,} if the following condition holds:
\[ c([\sigma]_i)\neq0,\forall\sigma\in\{1,2,\dots
,t\}^{\le m},  \ \textrm{and }\\ 1\leq i\leq m. \]
We note that, in particular, this condition implies that all of $c_1, \dots, c_t$ are non-zero
and hence $c_{\sigma} \neq 0$ for all $\sigma \in \{1, \dots, t\}^r$, for all $r \in \N$. Moreover,
when this condition is satisfied we let
\[\gamma := \gamma (c_1, \dots, c_t) :=  \min \{1, |c([\sigma]_i)|: \sigma \in \{1, \dots, t\}^{ \le m}, 1 \le i \le m\} >0.
\]

Below we collect some basic propositions which we will need for our main results. Throughout the rest
of this section we fix a sequence  $f_1, \dots, f_t$ of increasing functions which are pairwise almost disjoint.
We also fix $m \in \N$ such that $f_i(k) \neq f_j(k)$ for all $1 \le i < j \le t$ and $k \ge m$. And $c_1,\dots,c_t$
are some fixed constants in $\R$.

\begin{prop} \label{prp0} Let $\sigma, \tau \in\{1,2,\dots,t\}^r,$ $k, \ k' \in\N$ such that
 $f_{\sigma}(k)=f_{\tau}(k')$. Then $k=k'$.
\end{prop}
\begin{pf} We start with the case $r=1$. Let $\sigma=i$ and $\tau=j$. If $i=j$, then as $f_i$ is increasing, we have $k=k'$. If $i\neq j$, then $k=k'$ follows by the fact that $f_i$ and $f_j$ are almost disjoint. Suppose the proposition holds for $1,2,...,r$. Let $\sigma, \ \tau \in\{1,2,\dots,t\}^{r+1}$, $\sigma=i \sigma'$, $\tau=j \tau'$, with $\sigma', \ \tau' \in\{1,2,\dots,t\}^r,$ and assume $f_{\sigma}(k)=f_{\tau}(k')$.  So $ f_i(f_{\sigma'}(k))= f_{\sigma}(k)=f_{\tau}(k')=f_j(f_{\tau'}(k'))$. By step 1 of the induction we have $f_{\sigma'}(k)=f_{\tau'}(k')$. By step $r$ of the induction, we have $k=k'$.
\end{pf}
\begin{prop}\label{prp.25} Suppose $\sigma, \tau \in \{1,2,\cdots ,t\}^r$ and $k \ge m$. If $f_{\sigma}(k) = f_{\tau} (k)$,
then $\sigma = \tau$.
\end{prop}
\begin{pf}This simply follows from the fact that  $f_1(k), \dots, f_t(k) $ are disjoint for input values $k$
grater than $m$.
\end{pf}
\begin{prop}\label{prp0.5} Let $\sigma\in \{1,2,\dots ,t\}^r$, $r \ge m$, $ i \in \N$, and  $\sigma = \sigma_1 \sigma_2$
where $|\sigma_2| =m$. Then, $[\sigma]_i = \sigma_1 [\sigma_2]_i$.
\end{prop}
\begin{pf}
Let $\tau\in [\sigma]_i$ with $\tau=\tau_1\tau_2,$ and $|\tau_2| =m$. As $ f_{\sigma_1}(f_{\sigma_2}(i))=f_{\sigma}(i) = f_{\tau}(i) =f_{\tau_1}(f_{\tau_2}(i))$ and $|\sigma_1| = |\tau_1|$, by Proposition~\ref{prp0}, we have that $f_{\sigma_2}(i)=f_{\tau_2}(i).$ This means that $\sigma_2\sim_i\tau_2.$ As $f_{\sigma_2}(i) \ge m$  by Proposition~\ref{prp.25}, we have that $\sigma_1=\tau_1.$ Hence, we have shown that $\tau = \sigma_1 \tau_2$ where $\tau _2 \in [\sigma_2]$, completing the proof of the proposition.
\end{pf}

\begin{prop}\label{prp.75} Let $s> r \ge m$ and $i \in \N$. Then, $\sharp (s,i) =t^{s-r}\sharp (r,i)$.
\end{prop}
\begin{pf}This follows directly from Proposition~\ref{prp0.5}.
\end{pf}

\begin{prop} \label{prp2}
Let $\sigma\in\{1,2,\dots ,t\}^r$, $r\geq m.$ Let
$\sigma=\sigma_1\sigma_2,$ where
$\left|\sigma_2\right|=m.$ Then
$$c([\sigma]_i)=c_{\sigma_1}c([\sigma_2]_i).$$
\end{prop}
\begin{pf} This simply follows from Proposition~\ref{prp0.5}.
\end{pf}

\begin{prop} \label{prp3}
Assume that the coefficients $c_1,c_2, \dots , c_t$ satisfy the
non-zero condition at level $m$. Then, for all $\sigma\in\{1,2,\dots ,t\}^r, \textrm{with}\ r\geq m, \ \textrm{and }  i \in \N$ we have that
\[|c([\sigma]_i) | \ge  \gamma ^r\]
where $\gamma = \gamma (c_1, \dots, c_t)$. In particular, $|c([\sigma]_i) | >0$.
\end{prop}

\begin{pf} Let $\sigma\in\{1,2,\dots ,t\}^r, \textrm{with}\ r\geq m, \ \textrm{and} \ i\geq 1.$
Let us first consider the case that $i\geq m$.  Then, by Proposition~\ref{prp.25} we have that
$[\sigma]_i = \{\sigma\}$. Hence
$|c([\sigma]_i)|=|c_{\sigma}| \ge \gamma ^r $.

Now assume $1 \le i \le m$. Write $\sigma= \sigma_1 \sigma_2$ where $|\sigma_2| =m$.
Then, by Proposition~\ref{prp2}, we have that
\[ c([\sigma]_i) = c_{\sigma_1} c([\sigma_2]_i).
\]
Hence,
\[ |c([\sigma]_i) |= |c_{\sigma_1}|| c([\sigma_2]_i)| \ge \gamma ^{r-m} \cdot \gamma = \gamma ^{r-m+1} \ge \gamma ^r .\]

\end{pf}

\section{Main results}

We would like to prove the following result.

\begin{thm}\label{teore}
Let $f_1,f_2,\cdots,f_t$ be increasing functions which are
pairwise almost disjoint. Let $m\in \N$ be such that $f_i(k)\neq f_j(k)$ for all $1\leq i<j\leq t$ and $k\geq m.$
Assume that the coefficients $c_1,c_2,\cdots
c_t\in\R$ satisfy the non-zero condition at level $m$ and
let $T=\sum_{i=1}^{t}c_iT_{f_i}.$ Then, for sufficiently
large $\lambda$, the operator $\lambda T$ is chaotic.
\end{thm}

\begin{pf}
%[Proof.]

Let $\gamma = \gamma (c_1,\dots, c_t)$.
Let $\lambda$ be large enough so that
\begin{equation}\label{1}
\frac{4t^2}{\gamma}<\lambda.
\end{equation}
Let us first prove  that $\lambda T$ is transitive.
To this end, let $\vec{x}=(x_1,x_2, \dots, x_j, 0,0,0, \dots)$
 and $\vec{y}=(y_1,y_2, \dots, y_j, 0,0,0, \dots)$ be in  $l_p$ and let $\varepsilon >0 $  be given.
  It is enough to show that there exists
 $\vec{w}=(w_1,w_2, \dots )\in l_p$ such that $\|\vec{x}-\vec{w}\|_p ^p<\varepsilon$
 and $\|\vec{y}-\lambda^nT^n(\vec{w})\|_p^p <\varepsilon$ for some $n$.
 Without loss of generality we may assume that $j>m$.  Let $n$ be large enough so that  $n >j$ and so that
\begin{equation}\label{2}
j \cdot \|y\|_{\infty}^p\sum_{l=1}^{\infty}\frac{1}{4^{ln}} <\varepsilon.
\end{equation}
 Now we define $\vec{w}=(w_1,w_2, \dots )\in l_p$. Let $w_k=x_k,$ for $1\leq k\leq j.$ For each $l\geq 1,$  $\sigma\in\{1,2,\dots ,t\}^{ln}$ and $1\leq k\leq j,$ we define
\begin{equation}
\label{3} w_{f_{\sigma}(k)}=\frac{y_k\lambda^{-ln}}{c([\sigma]_k)\sharp (ln,k)}.
\end{equation}
 Let us observe that $c([\sigma]_k)\neq 0,$ by Proposition $\ref{prp3}.$  We let other $w_i$'s be zero.

 Let us first show that $\vec{w}$ is well defined.
Clearly, for $l\geq 1,$ $\sigma\in\{1,2,\dots , t\}^{ln}$, $1\leq k\leq j,$ we have that $f_{\sigma}(k)\geq ln\geq n>j.$ Hence, the definition of $w_{f_{\sigma}(k)}$ does not interfere with the definition of $\{w_1,w_2,\dots , w_j\}.$
 Suppose $f_{\sigma}(k)=f_{\sigma'}(k'),$ where $\sigma\in\{1,2,\dots ,t\}^{ln},$ $\sigma'\in\{1,2,\dots ,t\}^{l'n}$
 and $1\leq k,$ $k'\leq j,$ and $l, l'\geq 1.$ We need to verify that $w_{f_{\sigma}(k)}=w_{f_{\sigma'}(k')}.$ We first verify that $l =l'$. To obtain a contradiction, assume that $l \neq l'$.  Without loss of generality assume that $l<l'.$ Then, we write $\sigma'=\sigma_1'\sigma_2',$ where $|\sigma_1'|=|\sigma|.$ Note that $|\sigma_2'|\geq n$ and that $f_{\sigma}(k)=f_{\sigma'}(k')=f_{\sigma_1'}(f_{\sigma_2'}(k')).$ As $|\sigma|=|\sigma_1'|,$ we have that $k=f_{\sigma_2'}(k'),$ by Proposition \ref{prp0}. However, as $n>j$ we have that $ k= f_{\sigma_2'}(k')>j$ and this contradicts that $1\leq k \leq j.$ Hence $l=l'.$
 Now we have that $\sigma,\sigma'\in\{1,2,\dots , t\}^{ln}$ and since $f_{\sigma}(k)=f_{\sigma'}(k'), $ by Proposition \ref{prp0} we have that $k=k'.$ Hence $\sigma \sim_k \sigma'.$ Therefore,
 \[w_{f_{\sigma}(k)}=\frac{y_k\lambda^{-ln}}{c([\sigma]_k)\sharp (ln,k)}=\frac{y_{k'}\lambda^{-ln}}{c([\sigma']_{k'})\sharp (l'n,k')}  =w_{f_{\sigma'}(k')}.\]
 Therefore $\vec{w}$ is well defined.

Moreover,
\begin{eqnarray}
\label{eq:close1}\|\vec{x}\,-\,\vec{w}\|_p^p &= & \sum_{l=1}^{\infty}\sum_{k=1}^{j} \sum_{[\sigma ]_k  \in \A (ln, k)}
|w_{f_{\sigma}(k)}|^p   \\
\label{eq:close2}&=&\sum_{l=1}^{\infty}\sum_{k=1}^{j}\sum_{[\sigma ]_k  \in \A (ln, k)}\frac{|y_k|^p\lambda^{-pln}}{|c([\sigma]_k)|^p\sharp (ln,k)^p}\\
\label{eq:close3}&\le &\sum_{l=1}^{\infty}\sum_{k=1}^{j}\sum_{[\sigma ]_k  \in \A (ln, k)}\frac{\|y\|_{\infty}^p\lambda^{-pln} }{(\gamma^{ln})^p\sharp (ln,k)^p}\\
\label{eq:close4}&=&\sum_{l=1}^{\infty}\sum_{k=1}^{j}\frac{\|y\|_{\infty}^p\lambda^{-pln}  \sharp (ln,k)}{\gamma^{lnp}\sharp (ln,k)^p}\\
\label{eq:close5}&\le&j \|y \|_{\infty}^p \sum_{l=1}^{\infty} \frac{1}{(\gamma \lambda)^{pln}}\\
\label{eq:close6}&\le& j \|y \|_{\infty}^p \sum_{l=1}^{\infty} \frac{1}{4^{pln}} \le
j \|y \|_{\infty}^p \sum_{l=1}^{\infty} \frac{1}{4^{ln}}<\varepsilon.
\end{eqnarray}
Let us give explanations for some of the equalities and inequalities above. Let us first justify Equality~(\ref{eq:close1}).  Note that for all $i >j$ where $w_i \neq 0$ we have that $i = f_{\sigma}(k)$ for some $\sigma \in \{1, \dots, t\}^{nl}$, $l \ge 1$, and $1 \le k \le j$. Moreover, if $\sigma \sim_k \tau$
implies that $f_{\sigma}(k) = f_{\tau}(k)$. Hence, we only need to sum over one element from
each equivalence class in $\A(nl,k)$.
Equality~(\ref{eq:close2}) is simply inserting the definition of $w_{f_{\sigma}(k)}$. Inequality
~(\ref{eq:close3}) follows from Proposition~\ref{prp3}. Equality~(\ref{eq:close4}) follows from the definition of $\sharp(ln,k)$. Inequality~(\ref{eq:close5}) is simply algebra. Inequality~(\ref{eq:close6}) follows from the manner in which $\lambda$ and $n$ were chosen see (\ref{1}) and (\ref{2}) respectively.

Let us next show that $\|\lambda^nT^n(\vec{w})-\vec{y}\|_p^p<\varepsilon$. We first note that for all $k$
we have that
\[T^n(\vec{w})(k) = \sum_{\sigma \in \{1,2,\dots ,t\}^n} c_{\sigma}w_{f_{\sigma}(k)}.\]
In the case that $1 \le k \le j$, we have that
\begin{eqnarray}
 T^n(\vec{w})(k)& = &  \sum_{\sigma\in\{1,2,\dots ,t\}^{n}} c_{\sigma}w_{f_{\sigma}(k)}\nonumber\\
 & = & \sum_{\sigma\in\{1,2,\dots ,t\}^{n}} c_{\sigma}\frac{y_k\lambda^{-n}}{c([\sigma]_k)\sharp (n,k)}\nonumber\\
 &=&  \frac{y_k\lambda^{-n}}{\sharp (n,k)} \sum_{\sigma\in\{1,2,\dots ,t\}^{n}} \frac{ c_{\sigma}}{c([\sigma]_k)}\nonumber\\
 & = &  \frac{y_k\lambda^{-n}}{\sharp (n,k)} \sum_{[\sigma]_k \in \A(n,k)} \sum_{\tau \in [\sigma]_k} \frac{ c_{\tau}}{c([\tau]_k)}\nonumber\\
 \label{eq:nclose1}& =& \frac{y_k\lambda^{-n}}{\sharp (n,k)} \sum_{[\sigma]_k \in \A(n,k)} \sum_{\tau \in [\sigma]_k} \frac{ c_{\tau}}{c([\sigma]_k)}\\
 & =& \frac{y_k\lambda^{-n}}{\sharp (n,k)} \sum_{[\sigma]_k \in \A(n,k)} \frac{1}{c([\sigma]_k)}\sum_{\tau \in [\sigma]_k} { c_{\tau}}\nonumber\\
\label{eq:nclose2} &=& \frac{y_k\lambda^{-n}}{\sharp (n,k)} \sum_{[\sigma]_k \in \A(n,k)} \frac{1}{c([\sigma]_k)}\cdot {c([\sigma]_k)} = \frac{y_k\lambda^{-n}}{\sharp (n,k)} \sum_{[\sigma]_k \in \A(n,k)} 1 = y_k\lambda^{-n}.
 \end{eqnarray}
 We note that Equality~(\ref{eq:nclose1}) follows from the fact that $c([\sigma]_k) = c([\tau]_k)$ whenever $\sigma \sim_k \tau$. Equality~(\ref{eq:nclose2}) follows from the definitions of $c([\sigma]_k)$ and $\A(n,k)$.
 Hence we have just shown that for $1\le k \le j$, we have $\lambda^n T^n(\vec{w})(k) = y_k$.

 Let us now consider the case $k > j$. Let us compute $T^n(\vec{w})(k),$ for $k>j.$ In order to do that, we have to compute $w_{f_{\sigma}(k)}$ for $\sigma\in\{1,2,\dots ,t\}^n.$
If $w_{f_{\sigma}(k)}=0,\ \ \forall\sigma\in\{1,2,\cdots ,t\}^n,$ then $T^n(\vec{w})(k)=0.$
In the case in which there exists a $\sigma\in\{1,2,\cdots ,t\}^n$ such that  $w_{f_{\sigma}(k)}\neq0,$ we assign in a unique way to $k$ a triple $(l_k,\tau_k,i_k)$ such that $l_k\geq 1,$ $\tau_k\in\{1,2,\cdots ,t\}^{l_kn},$ $1\leq i_k\leq j$ and $k=f_{\tau_k}(i_k).$ Indeed as $w_{f_{\sigma}(k)}\neq 0,$ we have that $f_{\sigma}(k)=f_{\tau}(i_k), \ \ \textrm{for some}\ \ 1\leq i_k \leq j, \ \ \tau\in\{1,2,\dots , t\}^{ln}, \ \ l\geq 1.$
Now, let us write $\tau=\tau'\tau_k$, where $|\tau'|=n$ and $|\tau_k|=(l-1)n.$
Then, $f_{\sigma}(k)=f_{\tau'}(f_{\tau_k}(i_k)).$ Since $|\sigma|=|\tau'|$ we have that $k=f_{\tau_k}(i_k)$ by Proposition~\ref{prp0}.
As $k>j$ and $i_k\leq j,$ we have also that $\tau_k\neq\emptyset.$ Hence $l-1\geq 1.$ Let $l_k=l-1, $ therefore $(l_k,\tau_k,i_k)$ is the desired triple.

Now,
\begin{eqnarray}
|T^n(\vec{w})(k)|^p&=& \left|\sum_{\sigma\in\{1,2,\dots ,t\}^{n}} c_{\sigma}w_{f_{\sigma}(k)}\right|^p \nonumber\\
&\le &\sum_{\sigma\in\{1,2,\dots ,t\}^{n}}\left|c_{\sigma}w_{f_{\sigma}(f_{\tau_k}(i_k))}\right|^p\nonumber\\
&=&\sum_{\sigma\in\{1,2,\dots ,t\}^{n}} \left|c_{\sigma}w_{f_{\sigma\tau_k}(i_k)}\right|^p \nonumber\\
&=&\sum_{\sigma\in\{1,2,\dots ,t\}^{n}}\left| c_{\sigma} \frac{\lambda^{-(l_k+1)n}y_{i_k}}{c([\sigma\tau_k]_{i_k})\sharp((l_k+1)n;i_k)} \right|^p\nonumber\\
\label{eq:rclose1}&=& \sum_{\sigma\in\{1,2,\dots ,t\}^{n}}\frac{\lambda^{-l_knp}\lambda^{-np}|y_{i_k}|^p}{|c([\tau_k]_{i_k})|^p\sharp((l_k+1)n;i_k)^p}  \\
\label{eq:rclose2}&\le & \sum_{\sigma\in\{1,2,\dots ,t\}^{n}}\frac{\lambda^{-l_knp}\lambda^{-np}|y_{i_k}|^p}{\gamma^{l_knp}\sharp((l_k+1)n;i_k)^p}  \\
&\leq & \frac{t^n}{\lambda^{np}}\frac{\|y\|_{\infty}^p}{(\lambda \gamma)^{l_knp}}\nonumber
\end{eqnarray}
In above, Equality~(\ref{eq:rclose1}) follows by Proposition~\ref{prp2} and Inequality~(\ref{eq:rclose2}) follows by Proposition~\ref{prp3}.

Now let us make the following observation.
Suppose $k\neq k'>j,$  $T^n(\vec{w})(k)\neq 0$ and $T^n(\vec{w})(k')\neq 0.$ Then, $k\to (l_k,\tau_k,i_k)$ and $k'\to (l_{k'},\tau_{k'},i_{k'})$ are defined and, moreover, $(l_k,\tau_k,i_k)\neq (l_{k'},\tau_{k'},i_{k'}).$ Hence,

\begin{eqnarray}
\|\lambda^nT^n(\vec{w})-\vec{y}\|_p ^p& =&
\sum_{k=j+1}^{\infty} | \lambda^nT^n(\vec{w})(k)|^p= \lambda^{np}\sum_{k=j+1}^{\infty} | T^n(\vec{w})(k)|^p\nonumber\\
 &\leq &{\lambda^{np}}\sum_{l=1}^{\infty}\sum_{\tau\in\{1,2,\dots ,t\}^{ln}}\sum_{i=1}^j
\frac{t^n}{\lambda^{np}}\frac{\|y\|_{\infty}^p}{(\lambda \gamma)^{lnp}}\nonumber\\
&=& \sum_{l=1}^{\infty}t^{ln}\cdot j \cdot \frac{t^{n}\|y\|_{\infty}^p}{(\lambda \gamma)^{pln}} \nonumber\\
&\leq & j\|y\|_{\infty}^p\sum_{l=1}^{\infty}\frac{t^{(l+1)n}}{(\lambda \gamma)^{pln}}\nonumber \\
&\leq& j\|y\|_{\infty}^p \sum_{l=1}^{\infty}\frac{t^{2ln}}{(\lambda \gamma)^{pln}} \nonumber \\
\label{eq: basic1}&\le&  j \|y\|_{\infty}^p\sum_{l=1}^{\infty}
\left(\frac{t^2}{\lambda \gamma}\right)^{pln}<j\|y\|_{\infty}^p
\sum_{l=1}^{\infty}\left(\frac{1}{4}\right)^{pln}<\varepsilon.
\end{eqnarray}

In above, Inequality~(\ref{eq: basic1}) follows by (\ref{1}) and (\ref{2}).

Hence $\|\vec{y}-\lambda^nT^{n}(\vec{w})\|^p_p<\varepsilon$ and we have shown that $\lambda T$ is transitive.\\

Let us now prove that the periodic points of $\lambda T$ are dense in $\ell_p$. Indeed, $\vec{x}$ was chosen arbitrarily and $\|\vec{x} - \vec{w}\|_p^p < \varepsilon$. Hence, it will suffice to show that $\lambda^nT^n(\vec{w}) = \vec{w}$.

As we have already observed,  $\lambda^nT^n(\vec{w})(k)=y_k,\  \ \forall\ 1\leq k\leq j.$
Let us assume that $k>j.$ Let $A=\{f_{\sigma}(i):\sigma\in\{1,2,\dots ,t\}^{ln},\ 1\leq i\leq j,\ l\geq 1\}.$
 Note that the elements of $A$ are all greater than $j.$ Let $B=\N\setminus A\cap [j+1,\infty).$
Let us first consider the case $k\in A.$ Then $k=f_{\sigma}(i)$ for same $\sigma\in\{1,2,\dots , t\}^{ln}$, $1\leq i\leq j$ and $l\geq 1.$

 By definition we have that,
 \begin{eqnarray}
\lambda^nT^n(\vec{w})(k)
& = &  \lambda^n\sum_{\tau\in\{1,2,\dots ,t\}^{n}}c_{\tau}w_{f_{\tau}(k)}\nonumber\\
&=&\lambda^n\sum_{\tau\in\{1,2,\dots ,t\}^{n}}c_{\tau}w_{f_{\tau}(f_{\sigma}(i))}\nonumber\\
&=& \lambda^n\sum_{\tau\in\{1,2,\dots ,t\}^{n}}c_{\tau}w_{f_{\tau_{\sigma}}(i)}\nonumber\\
&=&\lambda^n\sum_{\tau\in\{1,2,\dots ,t\}^{n}}\frac{c_{\tau}\lambda^{-(l+1)n}y_i}{c([\tau\sigma]_i)\sharp ((l+1)n,i)}\nonumber\\
\label{eq:201}&=& \sum_{\tau\in\{1,2,\dots ,t\}^{n}}\frac{\lambda^{-ln}y_i}{c([\sigma]_i)\sharp ((l+1)n,i)}\\
\label{eq:202}&=&t^n\frac{\lambda^{-ln}y_i}{c([\sigma]_i)t^n\sharp (ln,i)}\\
&=&\frac{\lambda^{-ln}y_i}{c([\sigma]_i)\sharp (ln,i)}\nonumber\\
&=&w_k \nonumber.
\end{eqnarray}

We note that Equality~(\ref{eq:201}) follows from Proposition~\ref{prp2} and Equality~(\ref{eq:202}) follows
from Proposition~\ref{prp.75}.

Now let us assume that $k\in B.$ We claim that for each $\sigma\in\{1,2,\dots , t\}^{n}$ we have that $f_{\sigma}(k)\notin A.$
To obtain a contradiction, assume $f_{\sigma}(k)\in A$  and that there exist $l,\tau\in\{1,2,\dots , t\}^{ln}$ and $1\leq i\leq j$ such that $f_{\sigma}(k)=f_{\tau}(i).$ Let $\tau=\tau_1\tau_2,$ $|\tau_1|=n,$ then $f_{\sigma}(k)=f_{\tau_1}\left(f_{\tau_2}(i)\right).$ By the Proposition  $\ref{prp0}$ we have that $k=f_{\tau_2}(i).$
Moreover, as $k>j$, $\tau_2 \neq \emptyset.$ Hence $k\in A$ and this is a contradiction as $k\in B.$
Therefore, for all $k \in B$ we have that $f_{\sigma}(k) \in B$. Recall that the values of $\vec{w}$ at
each coordinate of $B$ is zero.
Hence, when $k \in B$ we have that
\[ \lambda^nT^n(\vec{w})(k)=\lambda^n\sum_{\sigma\in\{1,2,\dots ,t\}^{n}}c_{\sigma}w_{f_{\sigma}}(k)=0=w_k, \]
completing the proof of the theorem.
%\end{proof}

\begin{rem}
{\rm By the proof of the above Theorem it follows at once that if $f:\N\to\N$ is an increasing function, then for each $\lambda>1$ the operator $\lambda T_f$ is chaotic.}
\end{rem}

\begin{rem}
{\rm We observe that if the coefficients $c_1,c_2,\dots , c_t$ do not satisfy the non-zero condition at level m,  Theorem~\ref{teore} does not  hold. Indeed, let us consider the following increasing almost disjoint functions $f_1(n)=2n, \forall n\geq 1$ and $f_2(1)=2$ and $f_2(n)=2n-1, \  \forall n>1$.
The operator $T=2T_{f_1}-2T_{f_2}$ has the property that $T(\vec{x})$ is always zero in the first coordinate. Hence, $\lambda T$ is not chaotic for any $\lambda$.}
\end{rem}

\begin{cor} \label{c0}
Suppose that $f_1,f_2,\dots , f_t$ are increasing function with disjoint ranges. Let $c_1,c_2, \dots , c_t\in\R\setminus \{0\},$ then for sufficiently large $\lambda,$ the operator $\lambda\sum_{i=1}^tc_iT_{f_i}$ is chaotic.
\end{cor}

\begin{pf}
We note that  $c_1,c_2, \dots , c_t$ satisfies the non-zero condition at the level $m=1.$
The claim follows from Theorem~\ref{teore}.
\end{pf}

\begin{cor} \label{c1}
There exists  an infinite family ${\mathcal T}$ of chaotic operators such that for all $T_1,T_2, \dots , T_t\in{\mathcal T}$ and $c_1,c_2, \dots , c_t\in\R$ the operator $\lambda\sum_{i=1}^tc_iT_{f_i}$ is chaotic for sufficiently large $\lambda,$ provided that the following sum $$\sum_{i=1}^tc_iT_{f_i}$$ is not the zero operator.
\end{cor}
\begin{pf}

Let $f_1,f_2,\dots , f_i,\dots $ be increasing functions with disjoint ranges. Let ${\mathcal T}=\{2T_{f_i}:i\geq 1\}.$ Then ${\mathcal T}$ is the desired family.

\end{pf}
\begin{prop} \label{P0}
There exists a family ${\mathcal F}$ of cardinality continuum consisting of  increasing functions such that all distinct $f, g \in {\mathcal F}$ are pairwise almost disjoint.
\end{prop}

\begin{pf}
For each $c>1$ , let $g_c:\N \to \N$ be defined by $g_c(k)=\lceil ck \rceil$ where  $\lceil x \rceil$ denotes the ceiling function. Trivially $g_c$ is  increasing as $c>1.$ Moreover, if $c,d>1$ and $c\neq d$, then $\{k: g_c(k)=g_d(k)\}$ is finite.

Let $\alpha: \N^2\to\N$ be a bijection such that if
$(i,j), (i', j')\in\N^2$, with $i<i'$
and $j< j'$, then $\alpha (i,j) < \alpha (i', j')$.

For each $c>1$, define $f_c:\N\to\N$ by $f_c(k)=\alpha(k,g_c(k))$.

Let ${\mathcal F}=\{f_c:c>1\}$. First of all, let us note that $f_c$ is an  increasing map. Indeed, let $k\in\N$. Then, $f_c(k+1)=\alpha(k+1,g_c(k+1))>\alpha (k,g_c(k))=f_c(k)$.
Let $c,d\geq 1$, with $c\neq d$. As $\{k: g_c(k)=g_d(k)\}$ is finite and $\alpha$ is 1-1, we have that $\{k: f_c(k)=f_d(k)\}$ is finite. Let $k\neq k'$, then $(k,g_c(k))\neq (k', g_d(k'))$. As $\alpha$ is 1-1 we have that $f_c(k)=\alpha (k,g_c(k))\neq \alpha (k', g_d(k'))=f_d(k').$
Hence ${\mathcal F}$ has the desired proprieties.
\end{pf}

In the following, we use ${\mathcal L}^n$  to denote the $n$-dimensional Lebesgue measure in $\R^n$.
\begin{lem} \label{l0}

Let $P(x_1,x_2,\dots , x_t)$ be a non-zero polynomial in variables $x_1, \dots, x_t$. Then, the set $$\{(b_1,b_2,\dots , b_t)\in\R^t:P(b_1,b_2,\dots , b_t)=0\}$$ has $t$-dimensional Lebesgue measure zero.

\end{lem}

\begin{pf}
 We proceed by induction. The lemma is clear for $t=1$ since non-zero polynomials have only finitely many roots.

Suppose that the Lemma is true for $t=1, \dots, n.$  Now let us consider stage $t=n+1$. Let $P(x_1,x_2,\dots , x_{n+1})$ be a non-zero polynomial  in the variables $x_1, \dots, x_{n+1}$.
Fix $(r_1, \dots, r_{n+1}) \in \R^{n+1}$ such that
$P(r_1, \dots, r_{n+1}) \neq 0$.

Consider the polynomial $P(x_1,x_2,\dots , x_n,r_{n+1})$ in variables $x_1, \dots,x_n$.
Then this is a non-zero polynomial in $n$ variables and, by induction, we have that
\[A=\{(b_1,b_2,\dots , b_n)\in\R^n:P(b_1,b_2,\dots , b_n,r_{n+1})=0\}\]
has  $n$-dimensional Lebesgue measure zero.

 Now, let $B=\{ (b_1,b_2,\dots , b_n)\in\R^n:P(b_1,b_2,\dots , b_n, r_{n+1})\neq 0 \}$. For each $(b_1,b_2,\dots , b_n)\in B$, consider the polynomial of one variable $P(b_1, b_2,\dots, b_n, x_{n+1})$
 in variable $x_{n+1}$.  Then, this is a non-zero polynomial in one variable and, hence, the set $S_{b_1,b_2,\dots ,b_n} =\{y \in \R: P(b_1,b_2,\dots , b_n,y)=0\}$ is a finite set.

 Let us now make an observation. If $P(b_1,\dots ,  b_{n+1})=0$ and $P(b_1,\dots ,  b_{n}, r_{n+1}) =0$, then $(b_1,\dots ,  b_{n})  \in A$ and hence $(b_1,\dots ,  b_{n+1}) \subseteq A\times \R$.
 In the case, that $P(b_1,\dots ,  b_{n}, r_{n+1})  \neq 0$, we have that $(b_1,\dots ,  b_{n})  \in B$ and
 $b_{n+1} \in S_{b_1,b_2,\dots ,b_n}$. Putting all this together, we have that
  \begin{eqnarray*}
 \lefteqn{\{(b_1,\dots ,b_{n+1})\in\R^{n+1}:P(b_1,\dots ,  b_{n+1})=0\}}\\
 &&\subseteq (A\times \R)\cup\{(b_1,\dots , b_n, y)\in\R^{n+1}:(b_1,\dots , b_n)\in B, y\in S_{b_1,\dots ,b_n}\}.
\end{eqnarray*}
 As ${\mathcal L^n}(A) =0$, we have that \[{\mathcal L}^{n+1}(A\times \R)=0.\]
    By Fubini's Theorem, \[{\mathcal L}^{n+1}(\{(b_1,b_2,\dots , b_n, y): (b_1, \dots, b_n ) \in B, y\in S_{b_1,x_2,\dots ,b_n}\})=0.\]
 Hence, the set in question at the $(n+1)^{st}$ step of  induction has $(n+1)$-dimensional Lebesgue measure zero.\end{pf}

\begin{cor} \label{c2}
There exists a family ${\mathcal T}$ of cardinality continuum of chaotic operators such that for almost all $(c_1,c_2, \dots ,c_t)\in\R^t$ and $T_1, T_2, \dots , T_t\in \mathcal T$ the operator $\lambda\sum_{i=1}^tc_iT_{f_i}$ is chaotic for sufficiently large $\lambda.$
\end{cor}

\begin{pf}
Let ${\mathcal T}=\{2T_f:f\in{\mathcal F}\}$, where ${\mathcal F}$ is the family defined in the Proposition $\ref{P0}.$ Fix distinct elements in ${\mathcal F},$ let say $f_1, f_2, \dots , f_t.$ Choose $m$ according to the definition of almost disjoint functions and fix $\sigma\in\{1,2, \dots , t\}^{\leq m}$ and $1\leq i\leq j.$ Note that $c([\sigma]_i)$ is a non-zero polynomial in $c_1, c_2, \dots , c_t.$ By Lemma $\ref{l0}$, the $t$- dimensional Lebesgue measure of
$\{(c_1,c_2, \dots ,c_t): c([\sigma]_i)=0,\  \textrm{for some} \  \sigma\in \{1,2, \dots , t\}^{\leq m} \ \textrm{and}   1\leq i\leq m\}$ is equal to zero.
Hence by the Theorem $\ref{teore}$, for sufficiently large $\lambda$ and almost all $(c_1,c_2, \dots ,c_t)\in\R^t$, the operator $\lambda\sum_{i=1}^t2c_iT_{f_i}$ is chaotic.
\end{pf}

\begin{rem}
{\rm We observe that, since the vector space of finitely non-zero sequences is a dense subspace of $c_0$ (as well as of $l_p$, $1\leq p< \infty$) and since $l_p\subset c_0,$ then  Theorem $\ref{teore}$ and the other results of this section, hold also if the Banach space $X$ is $c_0$.}
\end{rem}

\end{pf}

\end{document}